\documentclass{article}
\usepackage{amsmath}
\usepackage{amssymb}
\usepackage{tikz-cd}
\usepackage[utf8]{inputenc}
\usepackage[english]{babel}
\usepackage[top=3cm, bottom=2cm, left=3cm, right=3cm]{geometry}
\usepackage{mathtools}
\usepackage{amsthm}
\usepackage{mathrsfs}
\usepackage{physics}
\usepackage{graphicx} 
\usepackage{esint}
\usepackage{hyperref}

\newtheorem{thm}{Theorem}

\newtheorem{dfn}[thm]{Definition}

\newenvironment{thmrecall}[1]
  {\par\medskip\noindent\textbf{Theorem~\ref{#1}.}\itshape}
  {\par\medskip}

\theoremstyle{Definition}

\title{On the size of the defect set of an analytic matrix field.}
\author{Aubertin Hugo}

\begin{document}
\maketitle
\begin{abstract}
    In this note, we study the defect set of a symmetric analytic matrix field $M$ defined on an open set $\Omega\subset \mathbb{R}^d$ with $d\ge2$. A defect point $x \in \Omega$ of $M$ is defined as a point for which no eigenframe can be chosen continuously in any neighborhood of $x$. We show that the defect set $S(M)$ is contained in a locally finite union of analytic submanifolds of codimension at least $2$. The proof follows fairly directly from results by Parusinski and Rond (Trans. Am. Math. Soc., 2020) and by Kurdyka and Paunescu (Duke Math. J., 2008).
    This result can be applied in the context of the Landau-de Gennes Theory to estimate the physical defects of a nematic liquid crystal.     
\end{abstract}

\section{Introduction}

A well-known result due to Rellich states that for an analytic and symmetric matrix field $M: ]-\epsilon, \epsilon[ \to Sym_n(\mathbb{R})$, there exists an analytic family of eigenvectors that diagonalizes $M$, see \cite{parusinski&rainer}. Our interest in the multiparameter setting is motivated by the physics of liquid crystals, where the study of the local behavior of the eigenvectors is closely related to physical defects. However, in this multiparameter setting, the study of the eigenframe becomes more complicated. Indeed, even for the basic example $M: \mathbb{R}^2 \to Sym_2(\mathbb{R})$ defined by \[
M(x,y)=\begin{pmatrix}
    x & y \\
    y & -x
\end{pmatrix},
\] there are no choices of continuous eigenvectors defined in a neighborhood of $0$ that diagonalize $M$. In this case, along the direction $y=0$, the canonical basis diagonalizes $M$ and along the direction $x=0$, the basis $ (\frac{1}{\sqrt{2}},\frac{1}{\sqrt{2}}),(\frac{1}{\sqrt{2}},-\frac{1}{\sqrt{2}})$ diagonalizes $M$. 

For a general matrix field $M: \Omega\subset \mathbb{R}^d \to Sym_n(\mathbb{R})$ that depends on more than one variable ($d\ge 2$), we define a defect point as a point $x \in \Omega$ for which no eigenframe can be chosen continuously in any neighborhood of $x$. We denote $S(M)$ the set of defect points of $M$. 

Although these obstructions prevent the existence of a continuous eigenframe in the original parameter space, they can be removed by replacing this space with a suitable modification $\sigma : W \to \Omega$, where $W$ is an analytic submanifold of dimension $d$ and $\sigma$ is locally an analytic diffeomorphism outside a nowhere dense set. More precisely, $\sigma$ is a locally finite composition of blowing-ups with smooth centers, see \cite{Paunescu&Kurdyka}.   

The first occurrence of this result comes from \cite{Paunescu&Kurdyka} where the authors show that we can locally construct an orthonormal eigenbasis that depends analytically on the parameter after performing two modifications.  

\medskip

More recently in \cite{parusinski&rond}, the authors show that only one modification is necessary. Indeed, they prove a diagonalization theorem for symmetric analytic matrix fields $M$ when the discriminant $\Delta_M$ is a monomial times a unit, that is of the form \[
\Delta_M(x)=x_1^{\alpha_1}\cdots x_d^{\alpha_d}u(x) 
\] with $u$ analytic that verify $u(0)\neq 0$. They then explain in their fifth section that after a modification, the discriminant can always be assumed to locally be of the appropriate form. 

The diagonalization theorem given in \cite{parusinski&rond} is closer to the question that motivates us : \begin{center}
    \textbf{What is the size of the defect set $S(M)$ ?}
\end{center}  Indeed, the theorem gives us a sufficient condition to be diagonalizable around a point, thus, the defect set will be contained in the set of points for which the discriminant of the characteristic polynomial of $M$ cannot be written in the form of a monomial times a unit in a local coordinate system.  

As mentioned earlier, this question comes from the physics of nematic liquid crystals, where the local orientation of rod-like molecules is modeled by a matrix field $Q:\Omega\to \mathcal{S}_0$ where $\mathcal{S}_0$ is the set of $3\times 3$ symmetric and traceless matrices. Equilibrium configuration minimize the Landau-de Gennes energy \[Q \mapsto \int_\Omega |\nabla Q|^2+f_b(Q)dx, \ Q \in H^1(\Omega, \mathcal{S}_0).\] Physical observations show the existence of points where the molecules fail to align in a well-defined direction. These points correspond to defects of the matrix field $Q$. Because of the physical modeling, modifying the space to remove the defects is not an option, this is why we focus on estimating the size of the defect set.

In the recent article \cite{Geng&Wang}, the authors consider minimizers of the Landau de Gennes energy under the Lyuksyutov constraint ($|Q|=1$) and $\Omega\subset \mathbb{R}^3$. They show that the subset of the  defects where the two first eigenvalues coincide, denoted $S_+(Q)$, is a $1$-rectifiable set, that is \[
S_+(Q)=\mathcal{C}\cup \mathcal{P},
\] where $\mathcal{C}$ is a countable union of $C^1$ curves and $\mathcal{P}$ has zero $\mathcal{H}^1$ Hausdorff measure.
To get this result, the authors study the first nonzero homogeneous term $ U$ in the Taylor expansion of $Q+\frac{\sqrt{6}}{2}(p\otimes p- \frac{1}{3}I_3)$ around a defect point $x\in S_+(Q)$, where $p$ is the third eigenvector of $Q$ and is well defined in a neighborhood of $x$. 

However, their method does not transfer to more general analytic and symmetric matrix fields, as in general, being defect-free is not stable under small perturbations. For example, consider the matrix field defined as \[
\tilde{M}(x,y)=M(x,y^2)=\begin{pmatrix}
    x & y^2 \\
    y^2 & -x
\end{pmatrix}.
\] The matrix field $\tilde{M}$ is the sum of a diagonal matrix $D=xI_2$ and a small perturbation $\tilde{M}-D=O(|(x,y)|^2)$. In this example, we see that the first nonzero term in the Taylor expansion $D$, has no defects, but that $\tilde{M}$ has a defect at $0$. Therefore, in general, we cannot show the absence of a defect by studying the first nonzero term of the Taylor expansion. 

\medskip

In this note, we show that arguments from \cite{Paunescu&Kurdyka} and \cite{parusinski&rond} imply rather directly a much stronger result on general analytic symmetric matrix fields. 
\begin{thm}\label{defect set of an analytic symmetric matrix field}
    Let $\Omega$ be an open subset of $\mathbb{R}^d$. Let $M:\Omega \to Sym_n(\mathbb{R})$ be a non zero analytic matrix field. Then, the defect set $S(M)$ is contained in a locally finite union of analytic submanifolds of codimension at least $2$.  
\end{thm}

In the context of liquid crystals, suitable boundary conditions are imposed. These conditions often confine the defects in a compact subset of $\Omega \subset \mathbb{R}^3$. In this case, the above theorem gives that the defects are contained in a finite union of analytic submanifolds of dimension at most $1$.

\medskip

The proof of the theorem is based on two main ingredients. The first one is from \cite{Paunescu&Kurdyka}, which implies that the discriminant is a monomial times a unit outside a closed semianalytic set of codimension at least $2$. The second one, as stated before, is the diagonalization theorem from \cite{parusinski&rond}.  

\medskip

To our knowledge, this result is new but its proof relies primarily on arguments already available in the literature. This is because the works from which these arguments come address different questions and do not directly consider the problem we are interested in.

\newpage
\section{Main result}
As discussed in the introduction, we assume that $d\ge2$ for the rest of this note, as there are no defects when $d=1$.

\begin{dfn}[{\cite[Def~2.1]{bierstone&milman}}]
Let $Z\subset \mathbb{R}^d$. We says that $Z$ is a semianalytic set if for any $x\in \mathbb{R}^d$, there exists some open neighborhood $U$ of $x$ together with a finite family of analytic functions $f_{ij}, g_{ij}: U\to \mathbb{R}$ such that \[
Z\cap U=\bigcup_{i=1}^{p}
\bigcap_{j=1}^{q_i}
\left\{
y \in U | \ 
f_{ij}(y)=0,\;
g_{ij}(y)>0
\right\}.
\]
Moreover, we say that $Z$ is an analytic set if we can choose $p=1$ and $g_{i,j}\equiv 1$. In other words, $Z$ is an analytic set if for any $x \in \mathbb{R}^d$, there exists some open neighborhood $U$ of $x$ together with a finite family of functions $f_i: U \to \mathbb{R}$ such that \[
Z\cap U = \{y \in U | \ f_1(y)=0,...,f_p(y)=0\}.
\]
\end{dfn}
The above definitions can be generalized to define semianalytic (and analytic) subsets of an analytic manifold $\mathcal{N}$ of dimension $d$.
The following theorem, due to Lojasiewicz, shows that any semianalytic set can be locally described by a finite union of pairwise disjoint analytic submanifolds.  

\begin{thm}[{\cite[Cor~2.11]{bierstone&milman}}]\label{Lojasiewicz Theorem}
    Let $Z\subset \mathbb{R}^d$ be a semianalytic set. Then, for any $x \in Z$, there exists an open neighborhood $U$ of $x$ such that \[
    Z\cap U =\bigcup_{i=0}^d \bigcup_{j=1}^{p_i}\Gamma_j^i,
    \]where $(\Gamma_j^i)_{j=1,...,p_i}$ are pairwise disjoint analytic submanifolds of dimension $i$. Moreover, we have that for any $1\le l\le d$ and $1\le m \le p_l$, \[
    \overline{\Gamma_m^l}\backslash \Gamma_m^l \subset \bigcup_{i=0}^{l-1} \bigcup_{j=1}^{p_i}\Gamma_j^i, 
    \] where $ \overline{\Gamma_m^l}$ denotes the closure of $ \Gamma_m^l$ in $U$. 
 \end{thm}

\begin{dfn}[{\cite[Rk~2.12]{bierstone&milman}}]
    Let $Z$ be a semianalytic set and $x \in Z$. We define the dimension of $Z$ at $x$, denoted $\dim_x(Z)$, by \[
    \dim_x(Z):=\max\{ \dim(\Gamma_j^i) | \ x \in \overline{\Gamma_j^i}\},
    \] 
    where $\Gamma_j^i$ comes from theorem \ref{Lojasiewicz Theorem}. We can define the dimension of $Z$ as the supremum of the dimensions of $Z$ at any of its points. The dimension at one point is well defined and does not depend on the choice of the stratification.

\end{dfn}

\begin{thm}[{\cite[Thm~4.2]{parusinski&rond}}]\label{diagonalization for monomial time unit discriminant}
    Let $\Omega  $ be an open subset of $\mathbb{R}^d$ and assume that $0\in \Omega$. Let $M: \Omega \to Sym_n(\mathbb{R})$ be an analytic matrix field. We assume that in a neighborhood of $0$, the discriminant $\Delta_M$ is a monomial times a unit, then there exists $\eta >0$ and $O:B(0,\eta) \to Sym_n(\mathbb{R})$ an analytic orthogonal matrix field such that for any $x \in B(0,\eta)$ \[
    O^{-1}(x)M(x)O(x)=\begin{pmatrix}
\lambda_1(x) & 0 & \cdots & 0 \\
0 & \lambda_2(x) & \ddots & \vdots \\
\vdots & \ddots & \ddots & 0 \\
0 & \cdots & 0 & \lambda_n(x)
\end{pmatrix},
    \] with $\lambda_i: B(0,\eta)\to \mathbb{R}$ analytic functions.     
\end{thm}

The above theorem is true if we replace $0$ with any point $a \in \Omega$ and $(x_1,...,x_d)$ with a local analytic coordinate system around $a$. 
We now recall the main result of this note. 

\begin{thmrecall}{defect set of an analytic symmetric matrix field}
    Let $\Omega$ be an open subset of $\mathbb{R}^d$. Let $M:\Omega \to Sym_n(\mathbb{R})$ be a non zero analytic matrix field. Then, the defect set $S(M)$ is contained in a locally finite union of analytic submanifolds of codimension at least $2$.    
\end{thmrecall}
 
\begin{proof}
    We may assume that $\Omega$ is connected. The proof is divided into two parts. In the first part, we show that the discriminant $\Delta_M: \Omega \to \mathbb{R}$ of the characteristic polynomial of $M$ is a monomial times a unit outside a semianalytic closed set $\Sigma \subset\Omega$ of codimension at least $2$. Then, we use the diagonalization theorem \ref{diagonalization for monomial time unit discriminant} to conclude that there are no defect points in $\Omega\backslash \Sigma$. The first part relies on an argument in \cite{Paunescu&Kurdyka}, more precisely from the proof of their theorem $5.5$. We reproduce the proof here for the reader's convenience.    

     We note that $\Delta_M$ is a nonzero analytic function, and thus its zero set $Z$ is an analytic set, we can moreover see that its dimension is at most $d-1$ using the identity principle. Let us denote by $Reg_{d-1}(Z)$ the set of regular points of dimension $d-1$ of $Z$, that is, the set of points $x\in Z$ such that $Z\cap B(x,\eta)$ is an analytic submanifold of dimension $d-1$ for $\eta>0$ small enough. By definition, it is clear that $Reg_{d-1}(Z)$ is open in $Z$. Then, the Lojasiewicz theorem \ref{Lojasiewicz Theorem} gives us that \[
    Z\backslash Reg_{d-1}(Z)\subset \Sigma'
    \] where $\Sigma'$ is a semianalytic closed set with dimension at most $d-2$. We now choose $C$, a connected component of $Reg_{d-1}(Z)$. Since $\Delta_M$ is nonzero and $\Omega$ is connected, we can choose $\alpha_C$ the smallest integer such that there exists a multi-index $(r_1,...,r_d)\in \mathbb{N}^d$ of length $\alpha_C$ with \[
     \left. \frac{\partial^{\alpha_C} \Delta(M) }{\partial_1^{r_1}\cdots \partial_d^{r_d}}  \right|_{C} \not \equiv 0.
    \]Let us denote $h_C:= \frac{\partial^{\alpha_C} \Delta(M) }{\partial_1^{r_1}\cdots \partial_d^{r_d}}$. Therefore, the set \[
    \Sigma''(C) := \{x \in \Omega | \ h_C(x)=0\}\cap C
    \] is a non trivial analytic set of $C$, and thus an analytic set of $\Omega$ with dimension at most $d-2$. Let $a \in C\backslash \Sigma''(C)$, to simplify the computation, we can assume that $a=0$ and $C= \{x_1=0\}$. Therefore, expressing $\Delta_M$ in the variables $x_1$ and $x'=(x_2,...,x_d)$, we have \[
    \Delta_M(x_1,x')=\sum_{k=0}^{+\infty}a_k(x')x_1^k.
    \] Hence, using the minimality of $\alpha_C$ and the fact that $ \Delta_M(0,\cdot)\equiv0$, we get that $\Delta_M(x_1,x')=x_1^{\alpha _C}u(x_1,x')$ for some analytic function $u$ that satisfies $u(0,0)\neq 0$. By denoting \[
    \Sigma'' = \bigcup \Sigma''(C)
    \]  the union over all the connected components $C$ of $Reg_{d-1}(Z)$, we see that \[
    \Sigma= \Sigma' \cup \Sigma''
    \] is a semianalytic closed set with codimension at least $2$, outside of which the discriminant $\Delta_M$ is a monomial times a unit in a local coordinate system.
    
    \medskip

    To conclude on the desired result, we see that any point $a \in \Omega \backslash \Sigma$ satisfies the hypothesis of theorem \ref{diagonalization for monomial time unit discriminant} in a local coordinate system, and therefore, we can choose analytically the eigenvectors and eigenvalues of $M$ in a small neighborhood of $a$. We get that \[
    S(M) \subset \Sigma,
    \] which is the desired result. 
\end{proof}

\newpage

\bibliographystyle{plain}
\bibliography{bibliographie.bib}

@article{Paunescu&Kurdyka,
 author = {Kurdyka, Krzysztof and Paunescu, Laurentiu},
 title = {Hyperbolic polynomials and multiparameter real-analytic perturbation theory},
 fjournal = {Duke Mathematical Journal},
 journal = {Duke Math. J.},
 issn = {0012-7094},
 volume = {141},
 number = {1},
 pages = {123--149},
 year = {2008},
 language = {English},
 doi = {10.1215/S0012-7094-08-14113-4},
 zbMATH = {5224875},
 Zbl = {1140.15006}
}

@article{parusinski&rond,
 author = {Parusi{\'n}ski, Adam and Rond, Guillaume},
 title = {Multiparameter perturbation theory of matrices and linear operators},
 fjournal = {Transactions of the American Mathematical Society},
 journal = {Trans. Am. Math. Soc.},
 issn = {0002-9947},
 volume = {373},
 number = {4},
 pages = {2933--2948},
 year = {2020},
 language = {English},
 doi = {10.1090/tran/8061},
 zbMATH = {7176218},
 Zbl = {1440.15009}
}

@article{Geng&Wang,
 author = {Geng, Zhiyuan and Wang, Changyou},
 title = {Eigenframe discontinuities of the {{\(Q\)}}-tensor model},
 fjournal = {Calculus of Variations and Partial Differential Equations},
 journal = {Calc. Var. Partial Differ. Equ.},
 issn = {0944-2669},
 volume = {64},
 number = {9},
 pages = {19},
 note = {Id/No 277},
 year = {2025},
 language = {English},
 doi = {10.1007/s00526-025-03158-4},
 zbMATH = {8105625},
 Zbl = {1578.35215}
}

@incollection{parusinski&rainer,
 author = {Parusi{\'n}ski, Adam and Rainer, Armin},
 title = {Perturbation theory of polynomials and linear operators},
 booktitle = {Handbook of geometry and topology of singularities VII},
 isbn = {978-3-031-68710-5; 978-3-031-68713-6; 978-3-031-68711-2},
 pages = {121--202},
 year = {2025},
 publisher = {Cham: Springer},
 language = {English},
 doi = {10.1007/978-3-031-68711-2_3},
 zbMATH = {8011697},
 Zbl = {1560.32028}
}

@article{bierstone&milman,
 author = {Bierstone, Edward and Milman, Pierre D.},
 title = {Semianalytic and subanalytic sets},
 fjournal = {Publications Math{\'e}matiques},
 journal = {Publ. Math., Inst. Hautes {\'E}tud. Sci.},
 issn = {0073-8301},
 volume = {67},
 pages = {5--42},
 year = {1988},
 language = {English},
 doi = {10.1007/BF02699126},
 url = {https://eudml.org/doc/104032},
 zbMATH = {4103444},
 Zbl = {0674.32002}
}
\end{document}